\documentclass{article}

\begin{document}

\title{Minimizing Squared Vertical and Squared Horizontal Errors}
\author{Donald E. Ramirez \and der@virginia.edu \and University of Virginia
\and Department of Mathematics \and P. O. Box 400137 \and Charlottesville,
VA 22904-4137}
\date{}
\maketitle

\begin{abstract}
The slope of the best fit line from minimizing the sum of both the squared
vertical errors and the squared horizontal errors is shown to be the root of
a fourth degree polynomial.
\end{abstract}

\section{Introduction}

With simple linear regression we have data $%
\{(x_{1},Y_{1}|X=x_{1}),...,(x_{n},Y_{n}|X=x_{n})\}$ and we minimize the sum
of the squared vertical errors. The question posed here is "Can we
effectively minimize both the sum of the squared vertical and squared
horizontal errors?" For notational convenience, we assume that the data is 
\textit{positively} correlated.

As an example, suppose we have paired data $(X,Y)$ where we first fit a
linear function $f(x)=y=\beta _{0}+\beta _{1}x$ to the data. For example, $Y$
could be the grade point average $GPA$ at graduation from a four year
university for a student, and $X$ is the corresponding $SAT$\ score before
matriculation$.$ Typically, Admissions Committees use such a least-squares
model to measure the effectiveness of the $SAT$ scores in the admission
process.

Suppose now we want to preform an inverse prediction at the value $y_{0}$ to
answer the question "What $SAT$ score should an admissions candidate receive
in order to have a predicted $GPA$ of, say, $2.0."$ This is found from the
inverse function $f^{-1}(y)=x=y/\beta _{1}-\beta _{0}/\beta _{1}.$

\section{Model}

For inverse prediction, we will want both $f(x)$ and $f^{-1}(y)$ to "fit"
the data, and we hope that the squared vertical and squared horizontal
errors will both be small for the fitted line $h(x)=\widehat{\beta }_{0}+%
\widehat{\beta }_{1}x$ which has minimized both the squared vertical and
squared horizontal errors. To that end, set 
\begin{equation}
SSE=\gamma \sum\limits_{i=1}^{n}(y_{i}-\beta _{0}-\beta
_{1}x_{i})^{2}+(1-\gamma )\sum\limits_{i=1}^{n}(x_{i}-y_{i}/\beta _{1}+\beta
_{0}/\beta _{1})^{2}  \label{SSE}
\end{equation}%
where $\gamma $ $(0\leq \gamma \leq 1)$. The parameter $\gamma $ allows for
a weighting of the two components of $SSE$ yielding the least square
estimators for $f(x)$ as $\gamma \rightarrow 1,$ and the least square
estimators of $f^{-1}(y)$ as $\gamma \rightarrow 0.$

We compute 
\begin{equation}
\frac{\partial }{\partial \beta _{0}}SSE=\frac{2\left( n\beta
_{0}-\sum\limits_{i=1}^{n}\left( y_{i}-\beta _{1}x_{i}\right) \right) \left(
\gamma \beta _{1}^{2}+1-\gamma \right) }{\beta _{1}^{2}}
\end{equation}%
with root%
\begin{equation}
\widehat{\beta }_{0}=\overline{y}-\widehat{\beta }_{1}\overline{x}
\label{b0 vetical and horizontal}
\end{equation}%
independent of $\gamma ,$ the same as in simple linear regression.

To find the slope $\hat{\beta}_{1},.$we compute%
\begin{eqnarray}
\frac{\partial }{\partial \beta _{1}}SSE &=&\gamma
\sum\limits_{i=1}^{n}\left( -2x_{i}y_{i}+2\beta _{0}x_{i}+2\beta
_{1}x_{i}^{2}\right)  \label{partial for b1} \\
&&+(1-\gamma )\left( \frac{-2n\beta _{0}^{2}}{\beta _{1}^{3}}%
+\sum\limits_{i=1}^{n}\left( \frac{2x_{i}y_{i}}{\beta _{1}^{2}}-\frac{2\beta
_{0}x_{i}}{\beta _{1}^{2}}-\frac{2y_{i}^{2}}{\beta _{1}^{3}}+\frac{4\beta
_{0}y_{i}}{\beta _{1}^{3}}\right) \right) .  \nonumber
\end{eqnarray}%
Set $S_{xx}=\sum_{i=1}^{n}(x_{i}-\overline{x})^{2},$ $S_{yy}=%
\sum_{i=1}^{n}(y_{i}-\overline{y})^{2}$ and $S_{xy}=\sum_{i=1}^{n}(x_{i}-%
\overline{x})(y_{i}-\overline{y})$, and let $\rho =S_{xy}/\sqrt{S_{xx}S_{yy}}
$ denote the correlation.

After some manipulation, the roots of Equation \ref{partial for b1} are
found by solving 
\begin{equation}
\gamma \sqrt{\frac{S_{xx}}{S_{yy}}}\beta _{1}^{4}-\gamma \rho \beta
_{1}^{3}+(1-\gamma )\rho \beta _{1}-(1-\gamma )\sqrt{\frac{S_{yy}}{S_{xx}}}%
=0.  \label{basic 4th degree}
\end{equation}

The (positive) root of Equation \ref{basic 4th degree} will be the slope of
the line which has minimized the $\gamma $-weighted sum of the squared
vertical and squared horizontal errors. 

With $\gamma =1.00$, the slope $\widehat{\beta }_{1}=\rho \sqrt{S_{yy}/S_{xx}%
};$ with $\gamma =0.00,$ the slope $\widehat{\beta }_{1}=(1/\rho )\sqrt{%
S_{yy}/S_{xx}}$; and in general, 
\begin{equation}
\rho \sqrt{S_{yy}/S_{xx}}\leq \widehat{\beta }_{1}\leq (1/\rho )\sqrt{%
S_{yy}/S_{xx}}  \label{bound}
\end{equation}

\section{An Example}

Suppose the data set is $\{(0,0),(0,0),(1,0),(1,1)\}$ with $\{\overline{x}%
=1/2,\overline{y}=1/4,S_{xx}=1,S_{yy}=3/4,\rho =\sqrt{3}/3=0.5774\}.$ If we
choose $\gamma =0.9,$ from Equation \ref{basic 4th degree}, $\widehat{\beta }%
_{1}=0.6612$; and from Equation \ref{b0 vetical and horizontal}, $\widehat{%
\beta }_{0}=1/4-(1/2)\widehat{\beta }_{1}=-.08060$. The bounds for $\widehat{%
\beta }_{1}$ are given in (\ref{bound}) and are $1/2\leq \widehat{\beta }%
_{1}\leq 3/2.$

\end{document}